\documentclass[aos,MSNbibl,seceqn,dvips]{arximspdf}
\usepackage{graphicx}
\doi{10.1214/13-AOS1197} 
\volume{42}
\issue{2}
\pubyear{2014}
\firstpage{654}
\lastpage{668}

\makeatletter
\def\longrightarrow{\rightarrow}
\newcommand{\rrVert}{\Vert}
\newcommand{\rrvert}{\vert}
\newcommand{\llVert}{\Vert}
\newcommand{\llvert}{\vert}
\def\limfunc{\operatorname}
\def\func{\operatorname}
\def\limfuncr{\mathrm}

\newtheorem{theorem}{Theorem}
\newtheorem{lemma}{Lemma}
\newtheorem{corollary}{Corollary}

\newproclaim{remark}{Remark}

\makeatother

\begin{document}
\begin{frontmatter}

\title{Oracally efficient estimation of autoregressive error distribution
with simultaneous confidence~band\thanksref{T1}}
\runtitle{Distribution estimation with simultaneous confidence band}
\thankstext{T1}{Supported in part by NSF award DMS-10-07594,
Jiangsu Specially-Appointed Professor Program SR10700111,
Jiangsu Province Key-Discipline (Statistics) Program ZY107002,
National Natural Science Foundation of China award 11371272,
Research Fund for the Doctoral Program of Higher Education
of China award 20133201110002 and Summer Fellowship of University of Toledo.}

\begin{aug}
\author[a]{\fnms{Jiangyan} \snm{Wang}\ead[label=e1]{wangjiangyan2007@126.com}},
\author[b]{\fnms{Rong} \snm{Liu}\ead[label=e2]{rong.liu@utoledo.edu}},
\author[c]{\fnms{Fuxia} \snm{Cheng}\ead[label=e3]{fcheng@ilstu.edu}}
\and
\author[d]{\fnms{Lijian} \snm{Yang}\corref{}\ead[label=e4]{yanglijian@suda.edu.cn}}
\affiliation{Soochow University, University of Toledo,
Illinois State University, and Soochow University and Michigan State
University}
\runauthor{Wang, Liu, Cheng and Yang}
\address[a]{J. Wang\\
Center for Advanced Statistics\\
\quad and Econometrics Research\\
Soochow University \\
Suzhou 215006\\
China \\
\printead{e1}}
\address[b]{R. Liu\\
Department of Mathematics\hspace*{16pt}\\
\quad and Statistics \\
University of Toledo\\
Toledo, Ohio 43606\\
USA\\
\printead{e2}}
\address[c]{F. Cheng\\
Department of Mathematics \\
Illinois State University \\
Normal, Illinois 61790\\
USA\\
\printead{e3}}
\address[d]{L. Yang\\
Center for Advanced Statistics\\
\quad and Econometrics Research \\
Soochow University \\
Suzhou 215006\\
China\\
and \\
Department of Statistics\\
\quad and Probability \\
Michigan State University \\
East Lansing, Michigan 48824\\
USA\\
\printead{e4}}
\end{aug}

\received{\smonth{11} \syear{2013}}
\revised{\smonth{12} \syear{2013}}

%
\begin{abstract}
We propose kernel estimator for the distribution function of unobserved
errors in autoregressive time series, based on residuals computed by
estimating the autoregressive coefficients with the Yule--Walker
method. Under mild assumptions, we establish oracle efficiency of the
proposed estimator, that is, it is asymptotically as efficient as the
kernel estimator of the distribution function based on the unobserved
error sequence itself. Applying the result of
Wang, Cheng and Yang
[\textit{J. Nonparametr. Stat.} \textbf{25} (2013) 395--407], the
proposed estimator is also asymptotically indistinguishable from the
empirical distribution function based on the unobserved errors. A
smooth simultaneous confidence band (SCB) is then constructed based on
the proposed smooth distribution estimator and Kolmogorov distribution.
Simulation examples support the asymptotic theory.
\end{abstract}

%


\end{frontmatter}

\section{Introduction}\label{sec1}

Consider an $\operatorname{AR} ( p ) $ process $\{X_{t} \}
_{t=-\infty}^{\infty}$ that satisfies
\[
X_{t}=\phi_{1}X_{t-1}+\cdots+
\phi_{p}X_{t-p}+Z_{t}
\]
in which $ \{ Z_{t} \} _{t=-\infty}^{\infty}$ are i.i.d. noises,
called errors, $\limfuncr{E}Z_{t}=0,\limfuncr{E}Z_{t}^{2}=\sigma^{2}$, with
probability density function (p.d.f.) $f ( z ) $ and cumulative
distribution function (c.d.f.) $F ( z ) =\int_{-\infty
}^{z}f (
u ) \,du$. For a positive integer $k$, the $k$-step ahead linear
predictor $\hat{X}_{n+k}$ of $X_{n+k}$, based on a length $n+p$
realization $%
\{ X_{t} \} _{t=1-p}^{n}$ up to time $n$, is well studied in
Chapters~5 and 9 of \cite{BD91}. While efficient methods are given to
compute $\hat{X}_{n+k}$ and its mean squared error, prediction
intervals are
unavailable unless the process is Gaussian; see Section~5.4 of~\cite{BD91}.

If $F ( z ) $ were known, all possible sample paths of the future
observation $X_{n+1}$ could be generated, and $P [ F^{-1} (
\alpha
_{1} ) \leq X_{n+1}-\hat{X}_{n+1}\leq F^{-1} ( \alpha
_{2} ) %
] =\alpha_{2}-\alpha_{1}$ for $0<\alpha_{1}<\alpha_{2}<1$. An
efficient estimator $\hat{F} ( z ) $ of $F ( z )
$ can be
used to construct a prediction interval $ [ \hat{X}_{n+1}+\hat{F}
^{-1} ( \alpha_{1} ),\hat{X}_{n+1}+\hat{F}^{-1} (
\alpha
_{2} )  ] $ for $X_{n+1}$, with confidence level $\alpha
_{2}-\alpha_{1}$. It is also pointed out in \cite{BN01} that
knowledge of
the c.d.f. $F ( z ) $ can improve related bootstrapping procedures.

While asymptotically normal estimators of the error density $f (
z ) $ have been studied in \cite{LN02,BD05} and \cite
{Cheng05a}%
, consistent estimator for error distribution $F ( z ) $
does not
exist for the $\operatorname{AR} ( p ) $ model. On the other hand, such
estimator has
been proposed for nonparametric regression in \cite{C02}, and
uniformly $%
\sqrt{n}$-consistent estimator of error distribution for the
nonparametric $\operatorname{AR}%
( 1 ) $--ARCH$ ( 1 ) $ model in \cite{NS13} and
nonparametric regression model in \cite{Cheng05b} and \cite{KN12}. It has
been used for symmetry testing in parametric nonlinear time series by~\cite%
{BN01}, and in nonparametric regression by~\cite{ND07}, as well as a
test of
parameter constancy in \cite{B96}. Other applications of error distribution
estimation include functional estimation: \cite{MSW04}; testing parametric
form of distribution and variance functions: \cite{NDN06} and \cite{DNV07};
testing for change-point in distribution: \cite{NV09} and testing for
additivity in regression: \cite{NV10} and \cite{MSW12}.

Assume for the sake of discussion that a sequence $ \{ Z_{t}
\}
_{t=1}^{n}$ of the errors were actually observed, \cite{Y73,F91,LY08} and
more recently \cite{WCY13} propose a kernel distribution
estimator (KDE) of $F ( z ) $ as
%
\begin{equation}
\tilde{F} ( z ) =\int_{-\infty}^{z}\tilde{f} ( u )
\,du=n^{-1}\sum_{t=1}^{n}\int
_{-\infty}^{z}K_{h} ( u-Z_{t} )
\,du,\qquad z\in\mathbb{R} \label{DEF:Ftilde}
\end{equation}
in which $K$ is a kernel function, with $K_{h} ( u ) =$ $%
h^{-1}K ( u/h ) $, and $h=h_{n}>0$ is called bandwidth. It
has been
established in \cite{F91} for Lipschitz continuous $F$, and for H\"{o}lder
continuous $F$ in \cite{WCY13} that $\tilde{F} ( z ) $ is
uniformly close to the empirical c.d.f. $F_{n} ( z ) $ at a
rate of $%
o_{p} ( n^{-1/2} ) $, thus inheriting all asymptotic
properties of
the latter. The general kernel smoothing results based on empirical process
in \cite{V94} require that $F\in C^{ ( 2 ) } ( \mathbb
{R}%
) $, thus excluding distributions such as the double exponential
distribution in our simulation study.

Unfortunately, $\tilde{F} ( z ) $ is infeasible, as one observes
only $ \{ X_{t} \} _{t=1-p}^{n}$, not $ \{ Z_{t}
\}
_{t=1}^{n}$. Denote by $\hat{\bolds{\phi}}$ the Yule--Walker
estimator of $%
\bolds{\phi}= ( \phi_{1},\ldots,\phi_{p} ) ^{T}$, then
%
\begin{eqnarray}\label{DEF:PHIHAT}
\hat{\bolds{\phi}} &=&\hat{\bolds{\Gamma}}_{p}^{-1}\hat{\bolds{\gamma}}%
_{p},\qquad \hat{\bolds{
\Gamma}}_{p}= \bigl\{ \hat {\gamma} ( i-j ) \bigr\} _{i,j=1}^{p},\qquad
\hat{\bolds{\gamma}}_{p}= \bigl( \hat{\gamma%
}(1),\ldots,
\hat{\gamma}(p) \bigr) ^{T},
\nonumber
\\[-8pt]
\\[-8pt]
\nonumber
\qquad\hat{\gamma}(l) &=&n^{-1}\sum_{i=1-p}^{n-\llvert  l\rrvert
}X_{i}X_{i+l},\qquad l=0,
\pm1,\ldots,\pm p.
\nonumber
\end{eqnarray}
We propose to estimate $F ( z ) $ by a two-step plug-in estimator
%
\begin{equation}
\hat{F} ( z ) =\int_{-\infty}^{z}\hat{f} ( u )
\,du=n^{-1}\sum_{t=1}^{n}\int
_{-\infty}^{z}K_{h} ( u-\hat{Z}%
_{t} ) \,du,\qquad z\in\mathbb{R} \label{DEF:Fhat}
\end{equation}
in which residuals $\hat{Z}_{t}=X_{t}-\hat{\phi}_{1}X_{t-1}-\cdots
-\hat{\phi%
}_{p}X_{t-p},1\leq t\leq n$.

Denote the empirical c.d.f.'s based respectively on $\hat{Z}_{t}$ and
$Z_{t}$ as
%
\begin{equation}
\hat{F}_{n} ( z ) =n^{-1}\sum_{t=1}^{n}I
\{ \hat{Z}%
_{t}\leq z \},\qquad F_{n} ( z )
=n^{-1}\sum_{t=1}^{n}I \{
Z_{t}\leq z \}. \label{DEF:ecdf}
\end{equation}

While $\hat{F}_{n} ( z ) $ is used for estimating $F (
z ) $, for example, in \cite{B96,BN01,C02,MSW04,Cheng05b,NDN06,DNV07,ND07,NV09,NV10,KN12,MSW12,NS13}, it is consistently shown
to be
less efficient than $F_{n} ( z ) $, as one referee observes;
see also Section~\ref{sec4.1}. Our
unique innovation is proving that the smooth estimator $\hat{F} (
z ) $ based on residuals is asymptotically equivalent to, not less
efficient than, the smooth estimator $\tilde{F} ( z ) $
based on
errors. As the Associate Editor points out, this result depends
crucially on the independence of $Z_{t}$ with $X_{t-r}$ for $r\geq1$,
ensured by the causal representation of the $X_{t}$ (proof of Lemma
\ref{order of uniform sumkxt}). We
have also learned from a referee that our result is related to the
orthogonality between innovation density $f ( z ) $ and
coefficient parameter $\bolds{\phi}$; see, for example, \cite{H86}.

Oracle efficiency of $\hat{F} ( z ) $ has powerful
implications, as
simultaneous confidence band (SCB) can be constructed for $F (
z ) $
over the entire real line, a natural tool for statistical inference on the
global shape of $F ( z ) $, which does not exist in previous works.
Working with a smooth estimator based on residuals can be adopted to other
settings such as nonparametric regression/autoregression, additive
regression, functional autoregression (FAR), etc., the present paper thus
serves as a first step in this direction.

Denote the distance between distribution functions as
%
\begin{eqnarray}\qquad
d ( F_{1},F_{2} ) &=&\llVert F_{1}-F_{2}
\rrVert _{\infty
}=\sup_{z\in\mathbb{R}}\bigl\llvert F_{1}
( z ) -F_{2} ( z ) \bigr\rrvert, \label{DEF:maxdev}
\\
D_{n} ( F_{n} ) &=&d ( F_{n},F ),\qquad D_{n} ( \hat{F}%
) =d ( \hat{F},F ),\qquad  D_{n}
( \tilde {F} ) =d ( \tilde{F},F ). \label{DEF:dnFn}
\end{eqnarray}

According to \cite{WCY13}, $d ( F_{n},\tilde{F} )
=o_{p} (
n^{-1/2} ) $, while it is well known that
%
\begin{equation}
P \bigl\{ \sqrt{n}D_{n} ( F_{n} ) \leq Q \bigr\}
\longrightarrow L(Q) \qquad\mbox{as } n\longrightarrow\infty, \label{DEF:kolmogorov}
\end{equation}
where $L ( Q ) $ is the classic Kolmogorov distribution function,
defined as
%
\begin{equation}
L(Q)\equiv1-2\sum_{j=1}^{\infty} ( -1 )
^{j-1}\exp \bigl( -2j^{2}Q^{2} \bigr),\qquad Q>0.
\label{DEF:limDn}
\end{equation}

Table~\ref{Tab:K-S test} displays the percentiles of $D_{n} (
F_{n} ) $ ($n\geq50)$, $L^{-1} ( 1-\alpha ) /\sqrt{n}$,
critical values for the two-sided Kolmogorov--Smirnov test.

\begin{table}
\caption{Critical values of Kolmogorov--Smirnov test}
\label{Tab:K-S test}
\begin{tabular*}{\textwidth}{@{\extracolsep{\fill}}lcccc@{}}
\hline
$\bolds{n}$ & $\bolds{\alpha=0.01}$ & $\bolds{\alpha=0.05}$ &
$\bolds{\alpha=0.1}$ & \multicolumn{1}{c@{}}{$\bolds{\alpha=0.2}$} \\
\hline
$\geq50$ & $1.63/\sqrt{n}$ & $1.36/\sqrt{n}$ & $1.22/\sqrt{n}$ &
$1.07/%
\sqrt{n}$ \\
\hline
\end{tabular*}
\end{table}

Theorem~\ref{THM:uniform1} entails that $d ( \hat{F},\tilde
{F} )
=o_{p} ( n^{-1/2} ) $, which together with \cite{WCY13},
lead to $%
\llvert \sqrt{n} \{ D_{n} ( \hat{F} ) -D_{n} (
F_{n} )  \} \rrvert \leq\sqrt{n} \{ d (
\hat{F},%
\tilde{F} ) +d ( F_{n},\tilde{F} )  \}
=o_{p} (
1 ) $. Applying Slutzky's theorem produces a smooth asymptotic
SCB by
replacing $\sqrt{n}D_{n} ( F_{n} ) $ in (\ref
{DEF:kolmogorov}) with
$\sqrt{n}D_{n} ( \hat{F} ) $.

The rest of the paper is organized as follows. Main theoretical results on
uniform asymptotics are given in Section~\ref{sec:theory}. Data-driven
implementation of procedure is described in Section~\ref
{sec:implementation}%
, with simulation results presented in Section~\ref{sec:examples}. Technical
proofs are in the \hyperref[app]{Appendix} and the supplemental
article~\cite{WLCY13}.

\section{Asymptotic results}
\label{sec:theory}

In this section, we prove uniform closeness of estimators $\hat
{F} (
z ) $ and $\tilde{F} ( z ) $ under H\"{o}lder continuity
assumption on $F$. For integer $\nu\geq0$ and $\beta\in (
0,1 ]
$, denote by $C^{ ( \nu,\beta ) } ( \mathbb{R}
) $ the
space of functions whose $\nu$th derivative satisfies H\"{o}lder condition
of order $\beta$,
%
\begin{equation}
C^{ ( \nu,\beta ) } ( \mathbb{R} ) = \biggl\{ \phi\dvtx %
\mathbb{R
\rightarrow R}\Big\llvert \sup_{x,y\in\mathbb{R}}\frac
{\llvert
\phi^{ ( \nu ) } ( x ) -\phi^{ ( \nu
)
} ( y ) \rrvert }{\llvert  x-y\rrvert ^{\beta
}}<+\infty \biggr\}. \label{DEF:Lipcont}
\end{equation}

We list some basic assumptions, where it is assumed that $\beta\in
(
1/3,1 ] $.

\begin{longlist}[(C1)]
\item[(C1)] \textit{The cumulative distribution function} $F\in
C^{ (
1,\beta ) } ( \mathbb{R} ) $, $0<f (
z )
\leq C_{f},\break \forall z \in\mathbb{R}$, \textit{where
}$C_{f}$ \textit{is a positive constant}.

\item[(C2)] \textit{The process} $ \{ X_{t} \} _{t=-\infty
}^{\infty}$ \textit{is strictly stationary with} $ \{
Z_{t} \}
_{t=-\infty}^{\infty}$ $\sim\break \limfunc{IID} ( 0,\sigma
^{2} ) $. $%
\{ X_{t} \} _{t=-\infty}^{\infty}$ \textit{is causal,
that is}, $\inf_{\llvert  z\rrvert \leq1}\llvert 1-\phi_{1}z-\cdots
-\phi_{p}z^{p}\rrvert >0$.

\item[(C3)] \textit{The univariate kernel function} $K ( \mathbf
{\cdot}%
) $ \textit{is a symmetric probability density, supported on}
$ [-1,1 ]$ \textit{and} $K\in C^{ ( 2 ) } (
\mathbb{R})$.

\item[(C4)] \textit{As} $n\rightarrow\infty$, $n^{-3/8}\ll
h=h_{n}\ll
n^{- \{ 2 ( 1+\beta )  \} ^{-1}}$.

\item[(C5)] $\limfuncr{E}\llvert  Z_{t}\rrvert ^{6+3\eta
}<\infty$,
\textit{for some} $\eta\in ( 6/5,+\infty ) $.
\end{longlist}

Conditions (C2), (C5) are typical for time series. Conditions (C1), (C3),
(C4) are similar to those in \cite{WCY13}. In particular, condition
(C4) on
bandwidth $h$ is rather different from those for constructing SCB in
\cite%
{BR73}.

The infinite moving average expansion $X_{t}=\sum_{j=0}^{\infty
}\psi_{j}Z_{t-j},t\in\mathbb{Z}$ and equation (3.3.6) of \cite{BD91}
ensure that there exist $C_{\psi}>0,0<\rho_{\psi}<1$, such that $%
\llvert \psi_{j}\rrvert \leq C_{\psi}\rho_{\psi}^{j},j\in
\mathbb{%
N}$. In particular,
\[
\bigl\{ \limfuncr{E}\llvert X_{t}\rrvert ^{6+3\eta} \bigr\}
^{1/ ( 6+3\eta ) }\leq\sum_{j=0}^{\infty}\llvert
\psi_{j}\rrvert \bigl\{ \limfuncr{E}\llvert Z_{t-j}\rrvert
^{6+3\eta} \bigr\} ^{1/ ( 6+3\eta ) }<\infty.
\]

In addition, the infinite moving average expansion and \cite{PT85} ensure
that there exist positive constants $C_{\rho}$ and $\rho\in (
0,1 ) $ such that $\alpha ( k ) \leq C_{\rho}\rho
^{k}$ holds for all~$k$, where the $k$th order strong mixing coefficient of the
strictly stationary process $ \{ X_{s} \} _{s=-\infty
}^{\infty}$ is defined as
\[
\alpha ( k ) =\sup_{B\in\sigma \{ X_{s},s\leq
t \},C\in\sigma \{ X_{s},s\geq t+k \} }\bigl\llvert P ( B\cap C ) -P ( B ) P
( C ) \bigr\rrvert,\qquad k\geq1.
\]

Our first result concerns asymptotic uniform oracle efficiency of $\hat{F}$
given in~(\ref{DEF:Fhat}) over intervals that grow to infinity with sample
size.

\begin{theorem}
\label{THM:uniform} Under conditions \textup{(C1)--(C5)}, the oracle estimator
$\hat{F}%
( z ) $ is asymptotically as efficient as the infeasible estimator
$\tilde{F} ( z ) $ over $z\in{}[-a_{n},a_{n}]$ where the
sequence $a_{n}>0$, $a_{n}\rightarrow\infty$, $a_{n}\leq C_{1}n^{C_{2}}$
for some $C_{1},C_{2}>0$, that is, as $n\rightarrow\infty$, $%
\sup_{z\in [ -a_{n},a_{n} ] }\llvert \hat{F} (
z ) -\tilde{F} ( z ) \rrvert =o_{p} (
n^{-1/2} ) $%
.
\end{theorem}

The above oracle efficiency of $\hat{F}$ extends to entire $\mathbb{R}$
provided that the extreme value of $ \{ \llvert  X_{t}\rrvert
\} _{t=1}^{n}$ has a mild growth bound. Denote
%
\begin{equation}
M_{n}=\max \bigl( \llvert X_{1}\rrvert,\llvert
X_{2}\rrvert,\ldots,\llvert X_{n}\rrvert \bigr)
\label{DEF:Mn}
\end{equation}

\begin{longlist}[(C6)]
\item[(C6)] \textit{There exists some} $\gamma>0$, \textit{such
that} $M_{n}=O_{p} ( n^{\gamma} )$.
\end{longlist}

Condition (C6) is satisfied, for instance, if the innovations $ \{
Z_{t} \} _{t=-\infty}^{\infty}$ have exponential tails. Denote
by $Z$
a random variable with the same distribution as the $Z_{t}$'s, the following
assumptions are as in A.1 and B.3 of \cite{R86}.

\begin{longlist}[(C5$'$)]
\item[(C5$'$)] \textit{There exist constants} $\sigma_{z}>0$
\textit{and} $\lambda>0$, \textit{such that one of the following conditions
holds}:
(1) $0<\lambda\leq1$ \textit{and for constants} $C_{z}>0,c_{z}\in
\mathbb{R}$,
%
\begin{equation}
P \bigl( \llvert \sigma_{z}Z\rrvert >z \bigr) \sim
C_{z}z^{c_{z}}\exp \bigl( -z^{\lambda} \bigr),\qquad z
\rightarrow+\infty; \label{EQ:pdf1}
\end{equation}
(2) $\lambda>1$ \textit{and for constants}
$C_{z,+},C_{z,-}>0,c_{z,+},c_{z,-}\in\mathbb{R}$,
\begin{eqnarray*}
\sigma_{z}^{-1}f \bigl( \sigma_{z}^{-1}z
\bigr) &\sim &C_{z,+}z^{c_{z,+}}\exp \bigl( -z^{\lambda} \bigr)
,\qquad z\rightarrow +\infty,
\\
\sigma_{z}^{-1}f \bigl( \sigma_{z}^{-1}z
\bigr) &\sim&C_{z,-} ( -z ) ^{c_{z,-}}\exp \bigl( - ( -z )
^{\lambda} \bigr),\qquad z\rightarrow-\infty.
\end{eqnarray*}
\textit{For} $D ( z ) =\sigma_{z}^{-1}f ( \sigma
_{z}^{-1}z ) \exp ( \llvert  z\rrvert ^{\lambda
} ) $
\textit{and its derivative} $D^{\prime} ( z ) $,
\begin{eqnarray*}
&&D ( z ) \sim C_{z,+}z^{c_{z,+}},\qquad z\rightarrow+\infty,\qquad
D ( z ) \sim C_{z,-} ( -z ) ^{c_{z,-}},\qquad z\rightarrow -
\infty,
\\
&&\qquad \limsup_{\llvert  z\rrvert \rightarrow\infty}\bigl\llvert zD^{\prime} ( z ) /D ( z )
\bigr\rrvert <\infty.
\end{eqnarray*}
\end{longlist}

Clearly the exponential tail condition (C5$'$) implies condition (C5), while
the next lemma establishes that it also entails condition (C6).

\begin{lemma}
\label{LEM:growthbound}Conditions \textup{(C2)}, \textup{(C5$'$)} imply $M_{n}=O_{p} (
( \log n ) ^{1/\lambda} ) $.
\end{lemma}

The next Theorem~\ref{THM:uniform1} extends Theorem~\ref{THM:uniform} with
the additional condition (C6) in general, or (C5$'$) in particular. As pointed
out by the associate editor, future works may lead to weaker conditions than
(C5$'$) that ensure (C6), aided by more powerful extreme value results
than in
\cite{R86}. We conjecture that Theorem~\ref{THM:uniform1} holds for
functional autoregression model (FAR) as well.

\begin{theorem}
\label{THM:uniform1} Under conditions \textup{(C1)--(C6)}, the oracle estimator
$\hat{F%
} ( z ) $ is asymptotically as efficient as the infeasible
estimator $\tilde{F} ( z ) $ over $z\in\mathbb{R,}$ that
is, as $%
n\rightarrow\infty,d ( \hat{F},\tilde{F} ) =\sup_{z\in
\mathbb{R}}\llvert \hat{F} ( z ) -\tilde{F} (
z )
\rrvert =o_{p} ( n^{-1/2} ) $. Especially, the above holds
under conditions \textup{(C1)--(C4)}, \textup{(C5$'$)}.
\end{theorem}

By \cite{B68}, as $n\rightarrow\infty$, $\sqrt{n} ( F_{n} (
z ) -F ( z )  ) \rightarrow_{d}B ( F
( z )
) $, where $B$ denotes the Brownian bridge. It is established in
\cite%
{WCY13} that $d ( F_{n},\tilde{F} ) =o_{p} (
n^{-1/2} ) $,
and hence Theorem~\ref{THM:uniform1} provides that $d ( \hat{F}%
,F_{n} ) =o_{p} ( n^{-1/2} ) $, and the following.

\begin{corollary}
\label{COR:smooth simultaneous cb} Under the conditions of Theorem~\ref%
{THM:uniform1}, as $n\rightarrow\infty$,
\[
\sqrt{n} \bigl( \hat{F} ( z ) -F ( z ) \bigr) \rightarrow _{d}B
\bigl( F ( z ) \bigr).
\]
For any $\alpha\in ( 0,1 ) $, $\lim_{n\rightarrow\infty
}P \{ F ( z ) \in\hat{F} ( z ) \pm
L_{1-\alpha}/%
\sqrt{n},z\in\mathbb{R} \} =1-\alpha$, and a smooth SCB for
$F (
z ) $ is
%
\begin{equation}
\bigl[ \max \bigl( 0,\hat{F} ( z ) -L_{1-\alpha}/\sqrt {n} \bigr),\min
\bigl( 1,\hat{F} ( z ) +L_{1-\alpha}/\sqrt {n} \bigr) \bigr],\qquad z\in\mathbb{R}
. \label{EQ:band}
\end{equation}
\end{corollary}

\section{Implementation}
\label{sec:implementation}

We now describe steps to construct the smooth SCB in~(\ref{EQ:band}).
For $%
n\geq50$, the following critical values from Table~\ref{Tab:K-S test} are
used:
\[
L_{1-0.01}=1.63,\qquad L_{1-0.05}=1.36,\qquad L_{1-0.1}=1.22,\qquad
L_{1-0.2}=1.07.
\]

To compute $\hat{F} ( z ) $ in (\ref{DEF:Fhat}), we use
the quartic
kernel $K(u)=15(1-u^{2})^{2}I \{ \llvert  u\rrvert \leq
1 \}
/16$ and a data-driven bandwidth $h=\func{IQR}\times n^{-1/3}$, with
$\func{%
IQR}$ denoting the sample inter-quartile range of $ \{ \hat{Z}%
_{t} \} _{t=1}^{n}$. This bandwidth satisfies condition (C4) as
long as
the H\"{o}lder order $\beta>1/2$. It is also similar to the robust and
simple one in \cite{WCY13}.

\section{Simulation examples}
\label{sec:examples}

In this section, we compare the performance of the estimator $\hat{F}$ with
the benchmark infeasible estimator $\tilde{F}$.

For sample sizes $n=50,100,500,1000$, a total of $1000$ samples $
\{
Z_{t} \} _{t=1}^{n}$ are generated, from the standard normal
distribution and the standard double exponential distribution, both of which
$\in C^{ ( 1,1 ) } ( \mathbb{R} ) $,
\[
F ( z ) =\int_{-\infty}^{z} ( 2\pi )
^{-1/2}e^{-u^{2}/2}\,du \quad\mbox{or}\quad F ( z ) =
\cases{1-1/2\exp ( -z ),&\quad $z\geq0,$\vspace*{2pt}
\cr
1/2\exp ( z ),&\quad $z<0$,}
\]
hence one would expect the data-driven bandwidth described in Section~\ref%
{sec:implementation} to perform well. We present results only for case 1:
standard normal distribution with the $\operatorname{AR} ( 1 ) $ model, and
case 2:
standard double exponential distribution with the $\operatorname{AR} ( 2 )
$ model.
Other combinations of error distributions and AR models have yielded similar
results which are omitted to save space.\looseness=-1

\subsection{Global errors}\label{sec4.1}

In this subsection, we examine the global errors of $\tilde{F}$ and~$\hat{F}$, measured by the maximal deviations $D_{n} ( \hat{F} )
,D_{n} ( \tilde{F} ) $ defined in (\ref{DEF:dnFn}), and
the Mean
Integrated Squared Error (MISE) defined as
%
\begin{eqnarray}
\limfunc{MISE} ( \hat{F} ) &=&\limfuncr{E}\int \bigl\{ \hat{F} ( z ) -F ( z )
\bigr\} ^{2}\,dz, \label{EQ:compute1}
\nonumber
\\[-8pt]
\\[-8pt]
\nonumber
\limfunc{MISE} ( \tilde{F} ) &=&\limfuncr{E}\int \bigl\{ \tilde{F}%
( z
) -F ( z ) \bigr\} ^{2}\,dz.
\end{eqnarray}

\begin{figure}

\includegraphics{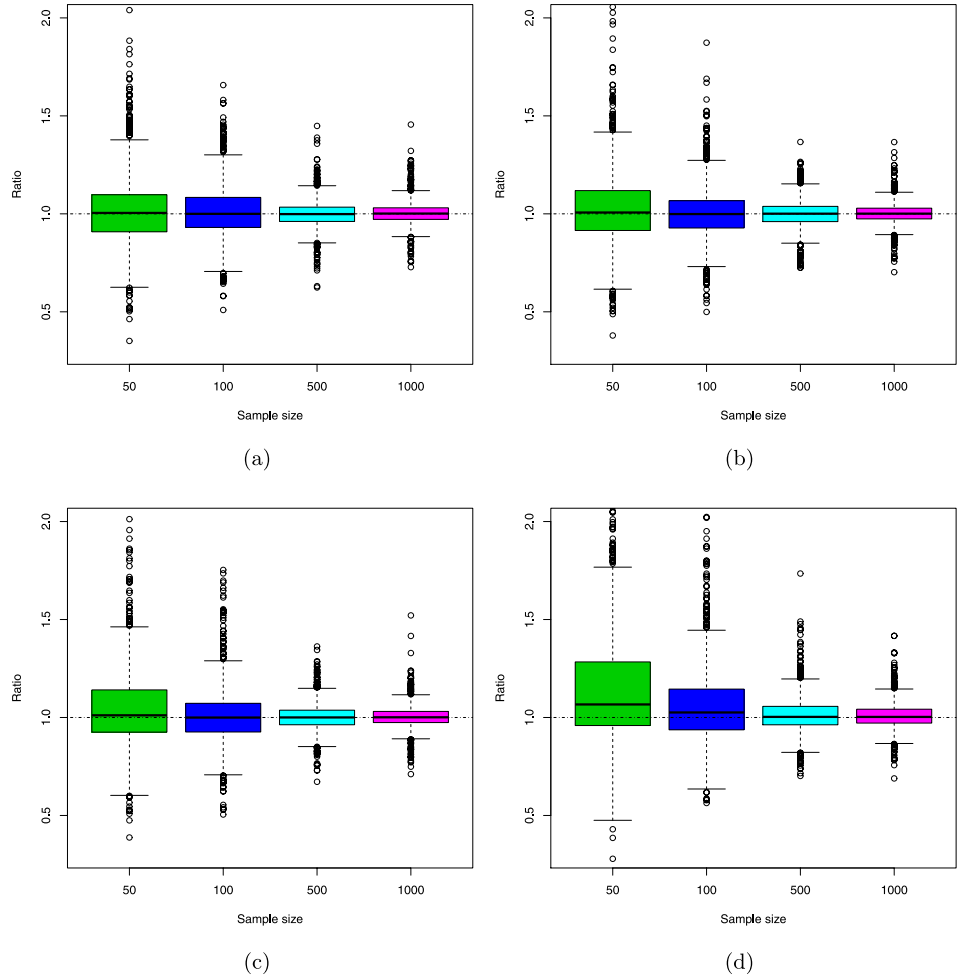}

\caption{Boxplot of the ratios $D_{n} ( \hat{F} )
/D_{n} (
\tilde{F} ) $ for $\operatorname{AR} ( 1 )$ model with standard normal
errors. The AR coefficients of \textup{(a)--(d)} are $-0.8$, $-0.2$, $0.2$, $0.8$,
respectively.}
\label{Fig: norm AR1 dn}
\end{figure}

Of interests are the means $\bar{D}_{n}(\hat{F})$ and $\bar
{D}_{n}(\tilde{F}%
) $ of $D_{n}(\hat{F})$ and $D_{n}(\tilde{F})$ over the $1000$ replications,
and similar means for $\limfunc{MISE}(\hat{F})$ and $\limfunc
{MISE}(\tilde{F}%
)$ for case~1. Table~\ref{Tab:norm AR1dnmse} contains these values, while
Figure~\ref{Fig: norm AR1 dn} is created based on the ratios
$D_{n}(\hat{F}%
)/D_{n}(\tilde{F})$ with four sets of coefficients. Both show that as $n$
increases, both deterministic ratios $\bar{D}_{n}(\hat{F})/\bar
{D}_{n}(%
\tilde{F})$ and $\overline{\limfunc{MISE}}(\hat{F})/\overline
{\limfunc{MISE}}%
(\tilde{F})\rightarrow1$, while the random ratio $D_{n}(\hat
{F})/D_{n}(%
\tilde{F})\rightarrow_{p}1$, all consistent with Theorem~\ref{THM:uniform1}.

\begin{table}
\caption{Comparing $\hat{F}$ and $\tilde{F}$: standard normal distribution
errors in $\operatorname{AR}(1)$}
\label{Tab:norm AR1dnmse}
\begin{tabular*}{\textwidth}{@{\extracolsep{\fill}}lccccc@{}}
\hline
$\bolds{\phi}$ & $\bolds{n} $ & $\bolds{\overline{D}_{n}(\hat{F}) }$ & $\bolds{\overline
{D}_{n}(\hat{F})/%
\overline{D}_{n}(\tilde{F})}$ & $\bolds{\limfunc{\overline{MISE}}(\hat
{F})}$ & $%
\bolds{\limfunc{\overline{MISE}}(\hat{F})/\limfunc{\overline
{MISE}}(\tilde{F})}$ \\
\hline
$-0.8$ & \phantom{00}$50 $ & $0.0857$ & $1.0012$ & $0.0028$ & $0.9881$ \\
& \phantom{0}$100 $ & $0.0649$ & $1.0034$ & $0.0015$ & $0.9927$ \\
& \phantom{0}$500 $ & $0.0306$ & $0.9966$ & $0.0003$ & $0.9983$ \\
& $1000 $ & $0.0228$ & $1.0022$ & $0.0002$ & $1.0012$ \\[3pt]
$-0.2$ & \phantom{00}$50 $ & $0.0865$ & $1.0105$ & $0.0029$ & $1.0213$ \\
& \phantom{0}$100 $ & $0.0646$ & $0.9985$ & $0.0015$ & $0.9927$ \\
& \phantom{0}$500 $ & $0.0307$ & $0.9997$ & $0.0003$ & $1.0020$ \\
& $1000 $ & $0.0228$ & $1.0027$ & $0.0002$ & $1.0048$ \\[3pt]
\phantom{$-$}$0.2 $ & \phantom{00}$50 $ & $0.0879$ & $1.0266$ & $0.0030$ & $1.0769$ \\
& \phantom{0}$100 $ & $0.0648$ & $1.0016$ & $0.0015$ & $1.0112$ \\
& \phantom{0}$500 $ & $0.0308$ & $1.0001$ & $0.0003$ & $1.0052$ \\
& $1000 $ & $0.0228$ & $1.0016$ & $0.0002$ & $1.0071$ \\[3pt]
\phantom{$-$}$0.8 $ & \phantom{00}$50 $ & $0.0977$ & $1.1418$ & $0.0046$ & $1.6205$ \\
& \phantom{0}$100 $ & $0.0688$ & $1.0631$ & $0.0019$ & $1.2780$ \\
& \phantom{0}$500 $ & $0.0311$ & $1.0122$ & $0.0003$ & $1.0615$ \\
& $1000 $ & $0.0230$ & $1.0094$ & $0.0002$ & $1.0333$ \\
\hline
\end{tabular*}
\end{table}

Table~4 in Section~2 of the supplemental article \cite{WLCY13}
contains $\bar{D}_{n}(\hat{F}), \linebreak\bar{D}%
_{n}(\hat{F}_{n}),\overline{\limfunc{MISE}}(\hat{F}),\overline{%
\limfunc{MISE}}(\hat{F}_{n})$ with $\hat{F}_{n}$ defined in (\ref
{DEF:ecdf}%
). Clearly, $\hat{F}$ outperforms $\hat{F}_{n}$ as we have commented
on page~3.

\subsection{Smooth SCBs}

In this subsection, we compare the SCBs based on smooth~$\hat{F}$, and the
infeasible $\tilde{F}$ and $F_{n}$ for case\vadjust{\goodbreak} 2, and $1-\alpha
=0.99,0.95,0.90,\linebreak0.80$. Table~\ref{Tab:Doubleexp cover AR2}
contains the coverage frequencies over $1000$ replications of the SCBs. The
smooth SCB is always conservative, the infeasible one more than the
data-based one in all cases except a few. The nonsmooth SCB based on $F_{n}$
has coverage frequencies closest to the nominal levels.

\begin{table}
\tabcolsep=0pt
\caption{Coverage frequencies for $\operatorname{AR} ( 2 ) $ model with double
exponential errors: left of parentheses-$\hat{F}$; right of
parentheses-$\tilde{F}$; inside the parentheses-$F_{n}$}
\label{Tab:Doubleexp cover AR2}
\begin{tabular*}{\textwidth}{@{\extracolsep{\fill}}lccccc@{}}
\hline
$\bolds{\phi}$ & $\bolds{n}$ & $\bolds{\alpha=0.01}$ &
$\bolds{\alpha=0.05}$ & $\bolds{\alpha=0.1}$ &
$\bolds{\alpha
=0.2}$ \\
\hline
& \phantom{00}$50$ & $0.998\ (0.995)\ 1.000$ & $0.992\ (0.973)\ 0.991$ &
$0.976\ (0.929)\ 0.980$ & $%
0.944\ (0.858)\ 0.950$ \\
$(-0.8$, & \phantom{0}$100$ & $0.998\ (0.992)\ 0.997$ & $0.987\ (0.963)\ 0.990$ & $%
0.972\ (0.924)\ 0.977$ & $0.928\ (0.858)\ 0.936$ \\
\hspace*{5pt}$-0.4)$ & \phantom{0}$500$ & $1.000\ (0.997)\ 1.000$ & $0.987\ (0.965)\ 0.984$ & $%
0.969\ (0.927)\ 0.965$ & $0.917\ (0.830)\ 0.923$ \\
& $1000$ & $0.995\ (0.992)\ 0.995$ & $0.985\ (0.951)\ 0.982$ &
$0.954\ (0.904)\ 0.949$ &
$0.889\ (0.814)\ 0.901$ \\[5pt]
& \phantom{00}$50$ & $0.994\ (0.995)\ 1.000$ & $0.981\ (0.973)\ 0.991$ &
$0.956\ (0.929)\ 0.980$ & $%
0.928\ (0.858)\ 0.950$ \\
$(0.8$, & \phantom{0}$100$ & $0.997\ (0.992)\ 0.997$ & $0.983\ (0.963)\ 0.990$ & $%
0.961\ (0.924)\ 0.977$ & $0.923\ (0.858)\ 0.936$ \\
\hspace*{2pt}$-0.4)$ & \phantom{0}$500$ & $1.000\ (0.997)\ 1.000$ & $0.982\ (0.965)\ 0.984$ & $%
0.966\ (0.927)\ 0.965$ & $0.914\ (0.830)\ 0.923$ \\
& $1000$ & $0.995\ (0.992)\ 0.995$ & $0.981\ (0.951)\ 0.982$ &
$0.950\ (0.904)\ 0.949$ &
$0.888\ (0.814)\ 0.901$ \\[5pt]
& \phantom{00}$50$ & $0.993\ (0.995)\ 1.000$ & $0.984\ (0.973)\ 0.991$ &
$0.965\ (0.929)\ 0.980$ & $%
0.930\ (0.858)\ 0.950$ \\
$(0.2$, & \phantom{0}$100$ & $0.998\ (0.992)\ 0.997$ & $0.983\ (0.963)\ 0.990$ & $%
0.961\ (0.924)\ 0.977$ & $0.915\ (0.858)\ 0.936$ \\
\hspace*{2pt}$-0.1)$ & \phantom{0}$500$ & $0.999\ (0.997)\ 1.000$ & $0.986\ (0.965)\ 0.984$ & $%
0.966\ (0.927)\ 0.965$ & $0.914\ (0.830)\ 0.923$ \\
& $1000$ & $0.995\ (0.992)\ 0.995$ & $0.982\ (0.951)\ 0.982$ &
$0.952\ (0.904)\ 0.949$ &
$0.896\ (0.814)\ 0.901$ \\[5pt]
& \phantom{00}$50$ & $0.991\ (0.995)\ 1.000$ & $0.979\ (0.973)\ 0.991$ &
$0.958\ (0.929)\ 0.980$ & $%
0.919\ (0.858)\ 0.950$ \\
$(0.2$, & \phantom{0}$100$ & $0.997\ (0.992)\ 0.997$ & $0.978\ (0.963)\ 0.990$ & $%
0.951\ (0.924)\ 0.977$ & $0.915\ (0.858)\ 0.936$ \\
\hspace*{6pt}$0.1)$ &\phantom{0}$500$ & $0.999\ (0.997)\ 1.000$ & $0.985\ (0.965)\ 0.984$ & $%
0.964\ (0.927)\ 0.965$ & $0.910\ (0.830)\ 0.923$ \\
& $1000$ & $0.995\ (0.992)\ 0.995$ & $0.983\ (0.951)\ 0.982$ &
$0.951\ (0.904)\ 0.949$ &
$0.894\ (0.814)\ 0.901$ \\
 \hline
\end{tabular*}
\end{table}

\begin{figure}

\includegraphics[scale=0.98]{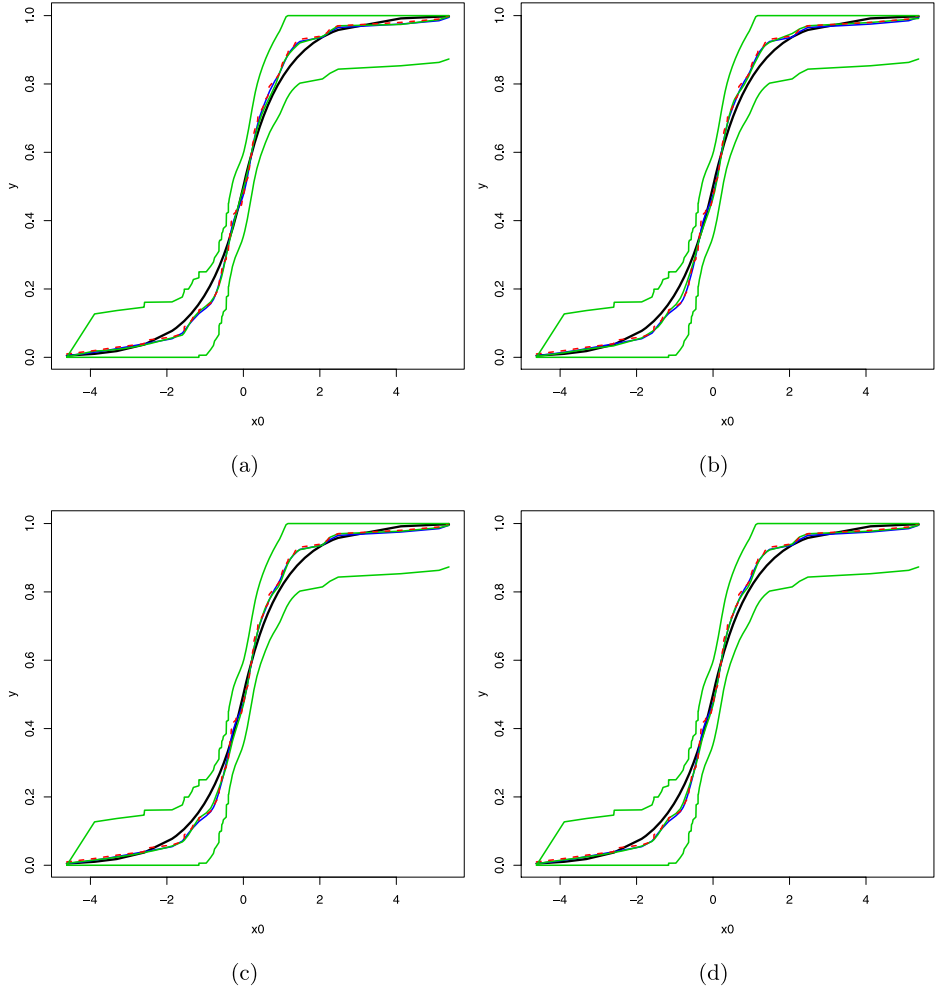}

\caption{Plots of the true c.d.f. $F$ (thick), the infeasible estimator
$\tilde{%
F}$ (solid), the data-based estimator $\hat{F}$ together with its
smooth $%
90\%$ SCB (solid) and $F_{n}$ (dashed) for $\operatorname{AR} ( 2 ) $ model
with $%
n=100$ standard double exponential errors. The AR coefficients of \textup{(a)--(d)}
are $ ( -0.8,-0.4 ) $, $ ( 0.8,-0.4 ) $, $ (
0.2,-0.1 ) $, $ ( 0.2,0.1 ) $, respectively.}
\label{Fig:doubleexp AR2 confidence bands}
\end{figure}

Figure~\ref{Fig:doubleexp AR2 confidence bands} depicts the true $F$
(thick), the infeasible $\tilde{F}$ (solid), the data-based $\hat{F}$ with
its $90\%$ SCB (solid) and $F_{n}$ (dashed), for a data of size $n=100$.
The three estimators are very close, with $\hat{F}$\ practically
distinguishable from $\tilde{F}$, consistent with our asymptotic theory.
Similar patterns have been observed for larger~$n$.

\begin{appendix}\label{app}
\section*{Appendix: Proofs}
\renewcommand{\thelemma}{A.\arabic{lemma}} %
\setcounter{equation}{0} 
\setcounter{lemma}{0}

\subsection{Preliminaries}

In this appendix, $C$ (or $c$) denote any positive constants, $U_{p}$
(or $%
u_{p}$) sequences of random variables uniformly $O$ (or $o$) of certain
order and by $O_{\mathrm{a.s.}}$ (or $o_{\mathrm{a.s.}}$) almost surely
$O$ (or $o$), etc.

The next two lemmas are used in the proof of Theorem~\ref{THM:uniform}.

\begin{lemma}[(\cite{B98}, Theorem~1.4)]\label{Bernstein} Let $\{\xi_{t}\}$ be a zero
mean real valued process. Suppose that there exists $c>0$ such that for
$%
i=1,\ldots,n$, $k\geq3$, $E\llvert \xi_{i}\rrvert ^{k}\leq
c^{k-2}k!E\xi_{i}^{2}<+\infty$, $m_{r}=\max_{1\leq i\leq N}\llVert \xi
_{i}\rrVert _{r},r\geq2$. Then for each $n>1$, integer $q\in
{}[
1,n/2]$, each $\varepsilon_{n}>0$ and $k\geq3$,
\begin{eqnarray*}
&& P \Biggl\{ \Biggl\llvert \sum_{i=1}^{n}
\xi_{i}\Biggr\rrvert >n\varepsilon _{n} \Biggr\}
\leq
a_{1}\exp \biggl( -\frac{q\varepsilon_{n}^{2}}{%
25m_{2}^{2}+5c\varepsilon_{n}} \biggr) +a_{2} ( k )
\alpha \biggl( %
\biggl[ \frac{n}{q+1} \biggr] \biggr)
^{{2k}/{(2k+1)}},
\end{eqnarray*}
where
$a_{1}=2\frac{n}{q}+2 ( 1+\frac{\varepsilon_{n}^{2}}{
25m_{2}^{2}+5c\varepsilon_{n}} ),a_{2} ( k )
=11n ( 1+%
\frac{5m_{k}^{2k/ ( 2k+1 ) }}{\varepsilon_{n}} )$.
\end{lemma}

\begin{lemma}[(\cite{BD91}, Theorem~8.1.1)]
\label{Phihat} The Yule--Walker
estimator $%
\hat{\bolds{\phi}}= ( \hat{\phi}_{1},\break \ldots, \hat{\phi}_{p} ) ^{T}$ of
$\bolds{\phi}=
( \phi_{1},\ldots,\phi_{p} ) ^{T}$ satisfies
$n^{1/2} ( \hat{\bolds{\phi}}-\phi ) \rightarrow N ( 0,\sigma^{2}\Gamma
_{p}^{-1} ) $, where $\Gamma_{p}$ is the covariance matrix
$ [
\gamma ( i-j )  ] _{i,j=1}^{p}$ with $\gamma (
h ) =%
\func{cov} ( X_{t},X_{t+h} ) $ for the causal $\operatorname{AR} (
p ) $
process $ \{ X_{t} \} $.\
\end{lemma}

\subsection{Proof of Theorem \texorpdfstring{\protect\ref{THM:uniform}}{1}}

\begin{lemma}
\label{order of alpha} Under conditions \textup{(C4)} and \textup{(C5)}, there exists an
$a>0$%
, such that the following are fulfilled for the sequence $ \{
D_{n} \} = \{ n^{a} \} $
%
\begin{eqnarray}\label{EQ:alpha}
\sum_{n=1}^{\infty}D_{n}^{-(2+\eta)}&<&
\infty,\qquad D_{n}^{-(1+\eta
)}n^{1/2}h^{1/2}
\rightarrow0,
\nonumber
\\[-8pt]
\\[-8pt]
\nonumber
 D_{n}n^{-1/2}h^{-1/2} ( \log n ) &\rightarrow&0.
\end{eqnarray}
\end{lemma}

\begin{lemma}
\label{order of uniform sumkxt} Under conditions \textup{(C1)--(C5)}, for any
$1\leq
r$, $s$, $v\leq p$,
%
\begin{equation}
\sup_{\llvert  z\rrvert \leq a_{n}}\Biggl\llvert n^{-1}\sum
_{t=1}^{n}K_{h} ( z-Z_{t} )
X_{t-r}\Biggr\rrvert =O_{\mathrm{a.s.}} \bigl( n^{-1/2}h^{-1/2}
\log n \bigr). \label{EQ:ukxt}
\end{equation}
\end{lemma}

\begin{lemma}
\label{order of uniform sumkprime xxt} Under conditions \textup{(C1)--(C5)},
for any
$1\leq r,s,v\leq p$,
%
\begin{eqnarray}
\sup_{\llvert  z\rrvert \leq a_{n}}\Biggl\llvert n^{-1}\sum
_{t=1}^{n}K_{h}^{\prime} (
z-Z_{t} ) X_{t-r}X_{t-s}\Biggr\rrvert
&=&O_{p} ( 1 ), \label
{EQ:kpxxt}
\\
\sup_{\llvert  z\rrvert \leq a_{n}}\Biggl\llvert n^{-1}\sum
_{t=1}^{n}K_{h}^{\prime\prime} (
z-Z_{t} ) X_{t-r}X_{t-s}X_{t-v}\Biggr
\rrvert &=&O_{p} ( 1 ). \label{EQ:kpxxxt}
\end{eqnarray}
\end{lemma}

\begin{lemma}
\label{order of uniform xxxxt} Under conditions \textup{(C1)--(C5)},
as $%
n\rightarrow\infty$,
\[
n^{-1}\sum_{t=1}^{n}\llvert
X_{t-r}X_{t-s}X_{t-v}X_{t-w}\rrvert
=O_{p} ( 1 ),\qquad 1\leq r,s,v,w\leq p.
\]
\end{lemma}

The proofs of Lemmas \ref{order of alpha}--\ref{order of uniform
xxxxt} are
in the supplemental article \cite{WLCY13}.

\begin{pf*}{Proof of Theorem~\ref{THM:uniform}}
Recall the definition of $\hat{Z}_{t}$ and $Z_{t}$ in
the \hyperref[sec1]{Intro-} \hyperref[sec1]{duction}; one has
%
\begin{eqnarray}\label{EQ:disf}
\hat{F} ( z ) -\tilde{F} ( z ) &=&n^{-1}\sum
_{t=1}^{n}\int_{ ( z-Z_{t} ) /h}^{ (
z-\hat{%
Z}_{t} ) /h}K
( v ) \,dv
\nonumber
\\[-8pt]
\\[-8pt]
\nonumber
&=& n^{-1}\sum_{t=1}^{n}
\biggl\{ G \biggl( \frac{z-\hat{Z}_{t}}{h} \biggr) -G \biggl( \frac
{z-Z_{t}}{h} \biggr)
\biggr\},
\end{eqnarray}
where $G ( z ) =\int_{-\infty}^{z}K ( u ) \,du$.
The right-hand side of equation (\ref{EQ:disf}) is
\begin{eqnarray*}
&&\frac{1}{n}\sum_{t=1}^{n} \biggl\{
G^{\prime} \biggl( \frac
{z-Z_{t}}{h%
} \biggr) \frac{\hat{Z}_{t}-Z_{t}}{h}+
\frac{1}{2}G^{\prime\prime
} \biggl( \frac{z-Z_{t}}{h} \biggr) \biggl(
\frac{\hat
{Z}_{t}-Z_{t}}{h} \biggr) ^{2}\\
&&\hspace*{95pt}{}+%
\frac{1}{6}G^{(3)}
\biggl( \frac{z-Z_{t}}{h} \biggr) \biggl( \frac
{\hat{Z}%
_{t}-Z_{t}}{h} \biggr)
^{3}+R_{t} \biggr\}.
\end{eqnarray*}

Therefore, $\hat{F} ( z ) -\tilde{F} ( z )
=I_{1}+I_{2}+I_{3}+I_{4}$, where
%
\begin{eqnarray}\label{EQ:parts}
I_{1} &=&n^{-1}\sum_{t=1}^{n}K
\biggl( \frac{z-Z_{t}}{h} \biggr) \frac{\hat{Z}_{t}-Z_{t}}{h},
\nonumber\\
I_{2} &=& ( 2n ) ^{-1}\sum_{t=1}^{n}K^{\prime}
\biggl( \frac{z-Z_{t}}{h} \biggr) \biggl( \frac{\hat
{Z}_{t}-Z_{t}}{h} \biggr)
^{2},
\\
I_{3} &=& ( 6n ) ^{-1}\sum_{t=1}^{n}K^{\prime\prime
}
\biggl( \frac{z-Z_{t}}{h} \biggr) \biggl( \frac{\hat
{Z}_{t}-Z_{t}}{h} \biggr)
^{3},\qquad I_{4}=n^{-1}\sum
_{t=1}^{n}R_{t}.
\nonumber
\end{eqnarray}

We now bound the four parts in (\ref{EQ:parts}).

Combining (\ref{EQ:parts}), Lemmas \ref{Phihat} and \ref{order of uniform
sumkxt}, for $1\leq r\leq p$,
%
\begin{eqnarray}\label{EQ:part1}
\sup_{\llvert  z\rrvert \leq a_{n}}\llvert I_{1}\rrvert &=&\sup_{\llvert  z\rrvert \leq a_{n}}n^{-1}\Biggl\llvert \sum
_{t=1}^{n}K \bigl\{ ( z-Z_{t} ) /h \bigr\}
\bigl\{ ( \hat{Z}_{t}-Z_{t} ) /h \bigr\} \Biggr\rrvert
\nonumber
\\[-8pt]
\\[-8pt]
\nonumber
&=&O_{p} \bigl( n^{-1/2} \bigr) O_{\mathrm{a.s.}} \bigl(
n^{-1/2}h^{-1/2}\log n \bigr) =o_{p} \bigl(
n^{-1/2} \bigr).
\end{eqnarray}

From (\ref{EQ:parts}), by applying Lemmas \ref{Phihat} and \ref
{order of
uniform sumkprime xxt}, for $1\leq r,s,v\leq p$,
%
\begin{equation}
\sup_{\llvert  z\rrvert \leq a_{n}}\llvert I_{2}\rrvert =o_{p}
\bigl( n^{-1/2} \bigr),\qquad \sup_{
\vert z\rrvert \leq a_{n}}\llvert
I_{3}\rrvert =o_{p} \bigl( n^{-1/2} \bigr).
\label{EQ:part2}
\end{equation}

According to (\ref{EQ:parts}), one has $\llvert  I_{4}\rrvert
\leq
n^{-1}\sum_{t=1}^{n}C\llvert  ( \hat{Z}_{t}-Z_{t} )
/h\rrvert ^{4}$. Thus,
%
\begin{eqnarray}\label{EQ:part4}
\sup_{z\in\mathbb{R}}\llvert I_{4}\rrvert &\leq& C\sup
_{1\leq t\leq
n}\bigl\llvert ( \hat{Z}_{t}-Z_{t} )
/h\bigr\rrvert ^{4}
\nonumber\\
&\leq&h^{-4}p^{4} \bigl( \max\llvert \phi_{r}-
\hat{\phi }_{r}\rrvert \bigr) ^{4}\sup n^{-1}\sum
_{t=1}^{n}\llvert X_{t-r}X_{t-s}X_{t-v}X_{t-w}
\rrvert
\\
&=&O_{p} \bigl( n^{-2}h^{-4} \bigr) \times
O_{p} ( 1 ) =o_{p} \bigl( n^{-1/2} \bigr),\qquad 1\leq
r,s,v,w\leq p,
\nonumber
\end{eqnarray}
which holds by using Lemmas \ref{Phihat} and \ref{order of uniform xxxxt}
simultaneously.

Since $\sup_{\llvert  z\rrvert \leq a_{n}}\llvert \hat{F}%
( z ) -\tilde{F} ( z ) \rrvert \leq
\sup_{_{\llvert  z\rrvert \leq a_{n}}} ( \llvert
I_{1}\rrvert +\llvert  I_{2}\rrvert +\llvert
I_{3}\rrvert
+\llvert  I_{4}\rrvert  ) $, Theorem~\ref{THM:uniform}
follows by
(\ref{EQ:part1}), (\ref{EQ:part2}) and (\ref{EQ:part4}) automatically.
\end{pf*}

\subsection{Proof of Lemma \texorpdfstring{\protect\ref{LEM:growthbound}}{1}}
Condition (C2) provides the infinite moving average expansion $%
X_{t}=\sum_{j=0}^{\infty}\psi_{j}Z_{t-j},t\in\mathbb{Z}$. Define
$\tilde{X}_{t}=\sum_{j=0}^{\infty}\llvert \psi_{j}\rrvert
\llvert  Z_{t-j}\rrvert $, so that $\llvert  X_{t}\rrvert \leq
\tilde{X}_{t},t\in\mathbb{Z}$. It is obvious that
%
\begin{equation}
M_{n}\leq\max ( \tilde{X}_{1},\tilde{X}_{2},
\ldots,\tilde{X}%
_{n} ). \label{EQ:max1}
\end{equation}

If $0$ $<\lambda\leq1$, then $C_{\psi}\rho_{\psi}^{j}=O (
j^{-\theta} ) $, for some $\theta>1$, $j\in\mathbb{N}$ and $%
\llvert \psi_{j}\rrvert \geq0$ according to condition (C5$'$), thus
condition A.1 of \cite{R86} is fulfilled, so Theorems~7.4 and~8.5 of
\cite{R86} imply that $\max ( \tilde{X}_{1},\tilde
{X}_{2},\ldots,%
\tilde{X}_{n} ) =O_{p} (  ( \log n ) ^{1/\lambda
} ) $%
. Thus, $M_{n}=O_{p} (  ( \log n ) ^{1/\lambda}
) $ by (%
\ref{EQ:max1}).

If $\lambda>1$, then $C_{\psi}\rho_{\psi}^{j}=O ( j^{-\theta
} ) $, for some $\theta>\max \{ 1,2 ( 1-1/\lambda
)
^{-1} \} $, $j\in\mathbb{N}$ according to condition (C5$'$), thus
condition B.3 of \cite{R86} is fulfilled. Theorem~6.1 in \cite{R86} implies
that $\max ( X_{1},\ldots,X_{n},-X_{1},\ldots,-X_{n} )
=O_{p} (  ( \log n ) ^{1/\lambda} ) $, hence $%
M_{n}=O_{p} (  ( \log n ) ^{1/\lambda} ) $.

Summarizing both scenarios, one concludes that under conditions (C2), (C5$'$),
$M_{n}=O_{p} (  ( \log n ) ^{1/\lambda} ) $, which
completes the proof.

\subsection{Proof of Theorem \texorpdfstring{\protect\ref{THM:uniform1}}{2}}

As in Theorem~\ref{THM:uniform}, equation (\ref{EQ:disf}) implies
\[
\sup_{_{z\in\mathbb{R}}}\bigl\llvert \hat{F} ( z ) -\tilde{F%
} (
z ) \bigr\rrvert \leq\sup_{_{z\in\mathbb{R}}} \bigl( \llvert
I_{1}\rrvert +\llvert I_{2}\rrvert +\llvert I_{3}
\rrvert +\llvert I_{4}\rrvert \bigr),
\]
where the four parts at the right-hand side are different from Theorem~\ref%
{THM:uniform} except $\sup_{_{z\in\mathbb{R}}}\llvert
I_{4}\rrvert $.\vspace*{1pt} So it remains to give the proof of parts $I_{1}$, $I_{2}$
and $I_{3}$ under conditions (C1)--(C6). The proof of next lemma is in the
supplemental article \cite{WLCY13}, where constants $\eta$ and
$\gamma$ are given in
conditions (C5) and (C6).

\begin{lemma}
\label{seckxt} Under conditions \textup{(C1)--(C6)}, for any $1\leq r,s,v\leq
p,a_{n}=h+n^{\delta}$, where $\delta> ( 7/4+6\gamma )
(
6+3\eta ) ^{-1}$
%
\begin{eqnarray}
\sup_{\llvert  z\rrvert >a_{n}}\Biggl\llvert n^{-1}\sum
_{t=1}^{n}K_{h} ( z-Z_{t} )
X_{t-r}\Biggr\rrvert& =&O_{p} \bigl( n^{-1} \bigr)
, \label{EQ:p1}
\\
\sup_{\llvert  z\rrvert >a_{n}}n^{-1}\Biggl\llvert \sum
_{t=1}^{n}K_{h}^{\prime} (
z-Z_{t} ) X_{t-r}X_{t-s}\Biggr\rrvert
&=&O_{p} \bigl( n^{-1} \bigr), \label{EQ:p2}
\\
\sup_{\llvert  z\rrvert >a_{n}}n^{-1}\Biggl\llvert \sum
_{t=1}^{n}K_{h}^{\prime\prime} (
z-Z_{t} ) X_{t-r}X_{t-s}X_{t-v}\Biggr
\rrvert &=&O_{p} \bigl( n^{-1} \bigr). \label{EQ:p3}
\end{eqnarray}
\end{lemma}

Theorem~\ref{THM:uniform1} is proved by combining Lemmas \ref{order of
uniform sumkxt}, \ref{order of uniform sumkprime xxt}, \ref{order of uniform
xxxxt} and \ref{seckxt}.
\end{appendix}

\section*{Acknowledgements}
This work is part of the first author's dissertation under the supervision
of the last author. The helpful comments from the Co-Editor, the Associate
Editor and two anonymous referees are gratefully acknowledged.


\begin{supplement}[id=suppA]
\stitle{Supplement to ``Oracally efficient estimation of autoregressive
error distribution with simultaneous confidence band''\\}
\slink[doi]{10.1214/13-AOS1197SUPP} 
\sdatatype{.pdf}
\sfilename{aos1197\_supp.pdf}
\sdescription{This supplement
contains additional technical proofs and some supporting numerical results.}
\end{supplement}

%

%
%
%
%
%

%

\printaddresses

\end{document}